\documentclass[notitlepage,leqno,10pt]{article}
\usepackage{latexsym}
\usepackage{amsmath}
\usepackage{amsfonts}
\usepackage{amssymb}

\newcommand{\intbar}{\mathop{\int\makebox(-13.5,0){\rule[4pt]{.7em}{0.3pt}}%
\kern-6pt}\nolimits}

\newcommand{\be}{\begin{equation}}
\newcommand{\ee}{\end{equation}}

\def\gt{\tilde{g}}

\def\sdtu{\sigma_2(\gt^{-1} A^1_{\gt})}
\def\sdu{\sigma_2(g^{-1} A^1_{g})}
\def\sutu{\sigma_1(\gt^{-1} A^1_{\gt})}
\def\suu{\sigma_1(g^{-1} A^1_{g})}
\def\sigut{\tilde{\sigma_1}}
\def\sigu{\sigma_1}
\def\sigdt{\tilde{\sigma_2}}
\def\sigd{\sigma_2}

\def\be{\begin{equation}}
\def\ee{\end{equation}}
\def\bea{\begin{eqnarray*}}
\def\eea{\end{eqnarray*}}
\def\f{\frac}

\author{Giovanni CATINO$^{a}$ and Zindine DJADLI$^{b}$}

\date{}

\title{Integral pinched $3$-manifolds are space forms}

\begin{document}

\newtheorem{lem}{Lemma}[section]
\newtheorem{pro}[lem]{Proposition}
\newtheorem{thm}[lem]{Theorem}
\newtheorem{rem}[lem]{Remark}
\newtheorem{cor}[lem]{Corollary}
\newtheorem{df}[lem]{Definition}
\newtheorem{claim}[lem]{Claim}
\newtheorem{conj}[lem]{Conjecture}

\maketitle

\begin{center}

%{\small Projet d'article - 26 janvier 2005}

\

{\small

\noindent $^a$ Universit\`a di Pisa - Dipartimento di Matematica

Largo Bruno Pontecorvo, 5

I56127 Pisa - Italy

\

\noindent $^b$ Institut Fourier - Universit\'e Grenoble 1

100 rue des Maths

F38402 Saint-Martin d'H\`eres Cedex - France

}

\end{center}

\footnotetext[1]{E-mail addresses: catino@mail.dm.unipi.it, \,
Zindine.Djadli@ujf-grenoble.fr}

\

\

\noindent {\sc abstract}. In this paper we prove that, under an
explicit integral pinching assumption between the $L^2$-norm of the
Ricci curvature and the $L^2$-norm of the scalar curvature, a closed
$3$-manifold with positive scalar curvature admits an Einstein
metric with positive curvature. In particular this implies that the
manifold is diffeomorphic to a quotient of ${\Bbb S}^3$.
\begin{center}

\bigskip\bigskip


\bigskip

\centerline{\bf AMS subject classification:  53C24, 53C20, 53C21, 53C25}

\end{center}

\centerline{}

\section{Introduction}\label{s:intro}

\noindent One of the basic questions concerning the relation between
algebraic properties of the curvature tensor and manifold topologies
is under which conditions on its curvature tensor a Riemannian
manifold is compact or homeomorphic to a space form (a manifold of
constant sectional curvature). For example, Bonnet-Myers theorem
states that a complete Riemannian manifold with positive lower bound
for its Ricci curvature is compact; the theorem of Klingenberg,
Berger and Rauch states that a compact, simply connected, $1\over
4$-pinched manifold with positive curvature is homeomorphic to the
standard sphere.

\medskip\noindent
In 1982, Hamilton \cite{Ha} introduced the Ricci flow and it appears
to be a very useful tool to study the relationships between topology
and curvature. For 2-dimensional compact manifolds, Hamilton
\cite{Ha3} and Chow \cite{Chow} proved that the normalized Ricci
flow converges and gave by the way a new proof of the well-known
uniformization theorem for compact surfaces. For 3 and 4-dimensional
compact manifolds with positive curvature, Hamilton, \cite{Ha} and
\cite{Ha2}, proved that the initial metric can be deformed into a
metric of constant positive curvature; it follows that these
manifolds are diffeomorphic to the sphere ${\Bbb S}^3$ or ${\Bbb
S}^4$, or a quotient space of ${\Bbb S}^3$ or ${\Bbb S}^4$ by a
group of fixed point free isometries in the standard metric. In
dimension 3, Hamilton's result is the following:

\begin{thm}[Hamilton]\label{Ha} If $(M,g)$ is a closed 3-dimensional Riemannian
manifold with positive Ricci curvature, then $M$ is diffeomorphic to
a spherical space form, i.e. $M$ admits a metric with constant
positive sectional curvature.
\end{thm}

\medskip\noindent
In this paper, we prove the existence of an Einstein metric of
positive curvature on compact, 3-dimensional manifolds satisfying an
integral pinching condition involving the second symmetric function
of the Schouten tensor.

\noindent More precisely, we consider $(M,g)$, a compact, smooth,
3-dimensional Riemannian manifold without boundary. Given a section
$A$ of the bundle of symmetric two tensors, we can use the metric to
raise an index and view $A$ as a tensor of type $(1,1)$, or
equivalently as a section of $End(TM)$. This allows us to define
$\sigma_2(g^{-1}A)$ the second elementary function of the
eigenvalues of $g^{-1}A$, namely, if we denote by $\lambda_1$,
$\lambda_2$ and $\lambda_3$ these eigenvalues
$$\sigma_2(g^{-1}A)=\lambda_1 \lambda_2+ \lambda_1\lambda_3
+\lambda_2\lambda_3.$$

\

\noindent In this paper we choose the tensor (here $t$ is a real number)
$$A^t_g=Ric_g-\f{t}{4}R_g g,$$
where $Ric_g$ and $R_g$ denote the Ricci and the scalar curvature of
$g$ respectively. Note that for $t=1$, $A^1_g$ is the classical
Schouten tensor $A^1_g=Ric_g-\f{1}{4}R_g g$ (see \cite{Besse}). Hence,
with our notations,
$\sigma_2(g^{-1} A^t_g)$ denotes the second elementary symmetric
function of the eigenvalues of $g^{-1}A^t_g$.

\

\noindent  Our present work is motivated by a recent paper of
M. Gursky and J. Viaclovsky \cite{GVbis}. Namely, they proved that,
giving a closed $3$-manifold $M$, a metric $g_0$ on $M$ (with normalized volume)
satisfying
$\int_M \sigma_2(g_0^{-1} A^1_{g_0})dV_{g_0} \ge 0$
is critical (over all metrics of normalized volume) for the functional
$${\mathcal F} \, : \, g \rightarrow \int_M \sigma_2(g^{-1} A^1_{g})dV_{g}$$
if and only if $g_0$ has constant sectional curvature.

\noindent Actually, it is not easy to exhibit a critical metric for this functional.
What we prove here (this is a consequence of
our main result in this paper) is that, assuming that there exists a metric $g$ on
$M$ with positive scalar curvature and such that
$\int_M \sigma_2(g^{-1} A^1_{g})dV_{g} \ge 0$ then the functional $\mathcal F$
admits a critical point (over all metrics of normalized volume)
$g_0$ with $\int_M \sigma_2(g_0^{-1} A^1_{g_0})dV_{g_0} \ge 0$.

\

\noindent We will denote $Y(M,[g])$ the Yamabe invariant associated to
$(M,g)$ (here $[g]$ is the conformal class of the metric $g$, that is
$[g]:=\left\{
{\tilde g}= e^{-2u}g \hbox{ for } u\in C^\infty(M)\right\}$). We recall that
$$Y(M,[g]) := \inf_{{\tilde g}\in [g]}
{\int_M R_{\tilde g} dV_{\tilde g} \over \left(\int_M dV_{\tilde g}\right)^{1\over 3}}.$$
An important fact that will be useful is that if $g$ has positive
scalar curvature then $Y(M,[g])>0$.

\medskip\noindent
Our main result is the following:
\begin{thm}\label{T2} Let $(M,g)$ be a closed 3-dimensional Riemannian
manifold with positive scalar curvature.
There exists a positive constant $C=C(M,g)$ depending only on
$(M,g)$ such that if
$$\int_M\sigd(g^{-1}A^1_g)\,dV_g+C\left({7\over 10} - t_0\right)Y(M,[g])^2 >0,$$
for some $t_0\leq 2/3$, then there exists a conformal metric
$\gt=e^{-2u}g$ with $R_{\gt}>0$ and
$\sigma_2(g^{-1}A^{t_0}_{\gt})>0$ pointwise. Moreover we have the
inequalities
\be\label{inequa}(3t_0-2)R_{\gt}\gt<6Ric_{\gt}<3(2-t_0)R_{\gt}\gt.\ee
\end{thm}

\noindent As an application, when $t_0=2/3$, we obtain
\begin{thm}\label{T1} Let $(M,g)$ be a closed 3-dimensional Riemannian
manifold with positive scalar curvature.
There exists a positive constant $C'=C'(M,g)$ depending only on
$(M,g)$ such that if
$$\int_M\sigd(g^{-1}A^1_g)\,dV_g+C'Y(M,[g])^2 >0,$$
then there exists a conformal metric $\gt=e^{-2u}g$ with positive
Ricci curvature ($Ric_{\gt} >0$). In particular if
$\int_M\sigd(g^{-1}A^1_g)\,dV_g \ge 0$ then there exists
a conformal metric $\gt=e^{-2u}g$ with positive Ricci curvature
($Ric_{\gt} >0$).
\end{thm}

\noindent Using Hamilton's theorem \ref{Ha}, we get:
\begin{cor}\label{corge}
 Let $(M,g)$ be a closed 3-dimensional Riemannian
manifold with positive scalar curvature.
There exists a positive constant $C'=C'(M,g)$ depending only on
$(M,g)$ such that if
$$\int_M\sigd(g^{-1}A^1_g)\,dV_g+C'Y(M,[g])^2 >0,$$
then $M$ is diffeomorphic to a spherical space form, i.e. $M$ admits
a metric with constant positive sectional curvature. In particular, if
$\int_M\sigd(g^{-1}A^1_g)\,dV_g \ge 0$ then $M$ is diffeomorphic
to a spherical space form.
\end{cor}

\begin{rem}\label{l:pinch}
Using the fact that $\sigd(g^{-1}A^1_g)=-{1\over 2}\left\vert
Ric_g\right\vert^2 + {3\over 16}R_g^2$, the assumption
$$\int_M\sigd(g^{-1}A^1_g)\,dV_g \ge 0$$
can be written
$$\int_M \left\vert Ric_g \right\vert^2 dV_g \le {3\over 8} \int_M R_g^2dV_g.$$
\end{rem}

\

\noindent Actually all these results are the consequence of the following
more general result:
\begin{thm}\label{Tgen} Let $(M,g)$ be a closed 3-dimensional Riemannian
manifold with positive scalar curvature.
There exists a positive constant $C=C(M,g)$ depending only on
$(M,g)$ such that if
$$\displaystyle \int_M\sigd(g^{-1}A^1_g)\,dV_g
+{1\over 24}\left({7\over 10} - t_0\right)
\inf_{g'=e^{-2u}g \, , \, \,  \vert \nabla_g u \vert_g \le C}\left(
\int_M R_{g'}^2e^{-u}dV_{g'}\right)>0,$$
for some $t_0\leq 2/3$, then there exists a conformal metric
$\gt=e^{-2u}g$ with $R_{\gt}>0$ and
$\sigma_2(g^{-1}A^{t_0}_{\gt})>0$ pointwise. Moreover we have the
inequalities
\be\label{inequal}(3t_0-2)R_{\gt}\gt<6Ric_{\gt}<3(2-t_0)R_{\gt}\gt.\ee
\end{thm}

\

\noindent There is a way to relate these result to the so-called $Q$-curvature
(the curvature associated to the Paneitz operator). The Paneitz operator introduced by
Paneitz in \cite{P} has demonstrated its importance in dimension 4 (see for example
Chang-Gursky-Yang \cite{cgyann} and \cite{cgy}). In dimension 3, the $Q$-curvature is defined by
$$Q_g = -{1\over 4}\Delta_g R_g - 2 \vert Ric_g \vert_g^2 +{23\over 32}R_g^2,$$
the Paneitz operator being  defined (in dimension 3) by
$$P_g = \Delta_g^2 - div_g \left(
-{5\over 4}R_g g + 4 Ric_g \right) d - {1\over 2}Q_g.$$
The Paneitz operator satisfies the conformal covariant property, that is, if
$\rho \in C^\infty(M)$, $\rho>0$, then for all $\varphi \in C^\infty(M)$,
$P_{\rho^{-4}g} (\varphi) = \rho^7 P_g(\rho\varphi)$. We can now state the Corollary:

\begin{cor}\label{Qcurv}
Let $(M,g)$ be a closed 3-dimensional Riemannian
manifold with non-negative Yamabe invariant. If there exists a metric
$g' \in [g]$ such that the $Q$-curvature of $g'$ satisfies
$$Q_{g'} \ge {1\over 48}R_{g'}^2,$$
then $M$ is diffeomorphic to a quotient of ${\Bbb R}^3$ if
$Y(M,[g])=0$ or to a spherical space form if $Y(M,[g])>0$.
\end{cor}
%According to this Corollary, if $M$ admits a metric $g'$
%such that the $Q$-curvature associated to $g'$ is non-negative
%and such that the Yamabe invariant of $M$ equipped with $g'$ is positive, then
%$M$ is diffeomorphic to a quotient of ${\Bbb S}^3$. This can be seen as a generalization to the
%dimension 3 of the uniformization of surfaces with  positive Gauss curvature.

\

\noindent Let us emphasize the fact that, in our results, we don't make any assumption on the
positivity of the Ricci tensor, we only assume that its trace is positive and a pinching on
its $L^2$-norm.

\noindent During the preparation of the manuscript
of this paper, we learned that Y. Ge, C.S. Lin and G. Wang
\cite{GLW} proved a weaker version of Corollary \ref{corge}, namely they prove that
if $(M,g)$ is a closed 3-dimensional Riemannian
manifold with positive scalar curvature
and if $\int_M\sigd(g^{-1}A^1_g)\,dV_g>0,$
then $M$ is diffeomorphic to a spherical space form. Their proof is completely
different from ours since they use a very specific conformal flow.

\

\noindent For the proof of Theorem \ref{T2} and Theorem \ref{T1},
we will be concerned
with the following equation for a conformal metric $\gt=e^{-2u}g$:
\be\label{eq1}
\left(\sigma_2(g^{-1}A^t_{\gt})\right)^{1/2}=fe^{2u},
\ee
where
$f$ is a positive function on $M$. Let $\sigma_1(g^{-1}A^1_g)$ be the trace of $A^1_g$
with respect to
the metric $g$. We have the following formula for the transformation
of $A_g^t$ under this conformal change of metric:
\be\label{trasf}
A^t_{\gt}=A^t_g+\nabla_g^2 u+(1-t)(\Delta_g u) g+du\otimes
du-\f{2-t}{2}|\nabla_g u|_g^2 g.
\ee
Since
$$A^t_g=A^1_g+(1-t)\sigma_1(g^{-1}A^1_g)g,$$
this formula follows easily from
the standard formula for the transformation of the Schouten tensor
(see \cite{Via}):
\be\label{trasf2}
A^1_{\gt}=A^1_g+\nabla_g^2
u+du\otimes du-\f{1}{2}|\nabla_g u|_g^2 g.
\ee
Using this formula we may
write (\ref{eq1}) with respect to the background metric $g$
$$\sigma_2\left(g^{-1}\left(A^t_g +\nabla_g^2 u+(1-t)(\Delta_g
u)g+du\otimes du-\f{2-t}{2}|\nabla_g u|_g^2
g\right)\right)^{1/2}=f(x)e^{2u}.$$

\

\

\noindent {\bf Aknowledgements : } The authors would like to
thank Sun-Yung Alice Chang and Paul Yang for their interest in
 this work.

\section{Ellipticity}
Following \cite{GV}, we will discuss the ellipticity properties of
equation (\ref{eq1}). \begin{df} Let
$(\lambda_1,\lambda_2,\lambda_3)\in {\Bbb R}^3$. We view the second
elementary symmetric function as a function on
${\Bbb R}^3$:
$$\sigma_2(\lambda_1,\lambda_2,\lambda_3)=\sum_{1\le i<j \le 3}\lambda_i\lambda_j,$$
and
we define
$$\Gamma^+_2=\{\sigma_2(\lambda_1,\lambda_2,\lambda_3)>0\}
\cap\{\sigma_1(\lambda_1,\lambda_2,\lambda_3)>0\} \subset {\Bbb R}^3,$$
where
$\sigma_1(\lambda_1,\lambda_2,\lambda_3)=\lambda_1+\lambda_2+\lambda_3$ denotes the trace.
\end{df}
For a symmetric linear transformation $A:V\rightarrow V$, where $V$
is an $n$-dimensional inner product space, the notation
$A\in\Gamma^+_2$ will mean that the eigenvalues of $A$ lie in the
corresponding set. We note that this notation also makes sense for a
symmetric 2-tensor on a Riemannian manifold. If $A\in\Gamma^+_2$, let
$\sigma_2^{1/2}(A)=\{\sigma_2(A)\}^{1/2}.$
\begin{df} Let $A:V\rightarrow V$, where $V$ is an $n$-dimensional inner product
space. The first Newton transformation associated with $A$ is (here $I$ is the identity map on $V$)
$$T_1(A):=\sigma_1(A)\cdot I-A.$$
Also, for $t\in {\Bbb R}$ we
define the linear transformation
$$\mathcal{L}^t(A):=T_1(A)+(1-t)\sigma_1(T_1(A))\cdot I.$$
\end{df}

\noindent We have the following:
\begin{lem}
If
$A \, : \, {\Bbb R}\rightarrow Hom(V,V)$, then
$$\f{d}{ds}\sigma_2(A)(s)=\sum_{i,j}T_1(A)_{ij}(s)\f{d}{ds}(A)_{ij}(s),$$
i.e, the first Newton transformation is what arises from
differentiation of $\sigma_2$.
\end{lem}
{\bf Proof} The proof of this lemma is a consequence of an easy computation.
See Gursky-Viaclovsky \cite{GVbis}

\begin{pro}(Ellipticity property)\label{Inver}
Let $u\in C^2(M)$ be a solution of equation
(\ref{eq1}) for some $t\leq 2/3$ and let
${\tilde g}=e^{-2u}g$. Assume that $A^t_{\gt}\in\Gamma^+_2$. Then
the linearized operator at $u$,
$\mathcal{L}^t:C^{2,\alpha}(M)\rightarrow C^{\alpha}(M)$, is
invertible $(0<\alpha<1)$.
\end{pro}
{\bf Proof} The proof of this proposition, adapted in dimension 3,
may be found in \cite{GV}.

\smallskip

\section{Upper bound and gradient estimate}

\noindent Throughout the sequel, $(M,g)$ will be a closed
3-dimensional Riemannian manifold with positive scalar curvature.
Since $R_g>0$, there exists $\delta>-\infty$ such that $A^\delta_g$ is
positive definite (i.e. $Ric_g - {\delta \over 4}R_g g >0$ on $M$).
Note that $\delta$ only depends on $(M,g)$.
For $t\in[\delta,2/3]$, consider the path of
equations (in the sequel we use the notation
$A^t_{u_t} := A^t_{g_t}$ for $g_t$ given by
$g_t=e^{-2u_t}g$)
\be\label{eq}
\sigma^{1/2}_2(g^{-1}A^t_{u_t})=fe^{2u_t},
\ee
where
$f=\sigma^{1/2}_2(g^{-1}A^\delta_g)>0.$ Note that $u\equiv 0$ is
a solution of (\ref{eq}) for $t=\delta$.

\begin{pro}[Upper bound]\label{UB} Let $u_t\in C^2(M)$ be a solution of
(\ref{eq}) for some $t\in[\delta,2/3]$. Then $u_t\leq\bar{\delta}$,
where $\bar{\delta}$ depends only on $(M,g)$.
\end{pro}
{\bf Proof} From Newton's inequality
$\sqrt{3}\sigma_2^{1/2}\leq\sigma_1$, so for all $x\in M$
$$\sqrt{3}fe^{2u_t}\leq\sigma_1(g^{-1}A^t_{u_t}).$$
Let $p\in M$ be a maximum of $u_t$, then using (\ref{trasf}), since
the gradient terms vanish at $p$ and
$(\Delta u_t)(p)\leq 0$,
 \bea
\sqrt{3}f(p)e^{2u_t(p)}&\leq&\sigma_1(g^{-1}A^t_{u_t})(p)\\&=&\sigma_1(g^{-1}A^t_{g})(p)+(4-3t)(\Delta
u_t)(p)\\&\leq&\sigma_1(g^{-1}A^t_{g})(p).
\eea
Since $t\geq\delta$,
this implies $u_t\leq\bar{\delta},$ for some $\bar{\delta}$ depending only on
$(M,g)$.

\begin{pro}[Gradient estimate]\label{GE} Let $u_t\in C^3(M)$ be a solution
of (\ref{eq}) for some $\delta\leq t\leq 2/3$. Assume that
$u_t\leq\bar{\delta}$. Then $\parallel\nabla_g
u\parallel_{{g,\infty}}<C_1$, where $C_1$ depends only on
$(M,g)$ and $\bar\delta$.
\end{pro}

\noindent The proof of this lemma can be found in the paper Gursky-Viaclovsky \cite{GV}.

\begin{rem}
Note that we will use this proposition with $\bar\delta$ given
by Proposition \ref{UB} and then, since $\bar\delta$ depends only
on $(M,g)$, we infer that $C_1$ only depends on $(M,g)$.
\end{rem}

\section{A technical lemma}

\noindent As we proved in the previous section, there exists two constants
$\bar\delta$ and $C_1$ depending only on $(M,g)$ such that all solutions of
(\ref{eq}) for some $\delta\leq t\leq 2/3$, satisfying
$u_t\leq\bar{\delta}$ satisfies $\parallel\nabla_g
u\parallel_{{g,\infty}}<C_1$.

\noindent We consider the following quantity:
$$I(M,g) :=
\inf_{g'=e^{-2\varphi}g \, , \, \,  \vert \nabla_g \varphi \vert \le C_1}\left(
\int_M R_{g'}^2e^{-\varphi}dV_{g'}\right).$$
We let, for $g'=e^{-2\varphi}g$
$$i(g'):=\int_M R_{g'}^2e^{-\varphi}dV_{g'}.$$
As one can easily check, if two metrics
$g_1$ and $g_2$ are homothetic, then $i(g_1)=i(g_2)$. So, we have
$$I(M,g) =
\inf_{g'=e^{-2\varphi}g \, , \, \,  Vol (M,g')=1 \,\, \hbox {\small and} \,\, \vert \nabla_g \varphi \vert_g \le C_1}\left(
\int_M R_{g'}^2e^{-\varphi}dV_{g'}\right).$$

\noindent We have the following
\begin{lem}\label{techn}
There exists a positive constant $C=C(M,g)$ depending only on
$(M,g)$ such that
$$I(M,g) \ge C \left( Y(M,[g])\right)^2.$$
\end{lem}

\noindent {\bf Proof} As we have seen
$$I(M,g) =
\inf_{g'=e^{-2\varphi}g \, , \, \,  Vol (M,g')=1 \,\, \hbox {\small and} \,\, \vert \nabla_g \varphi \vert_g \le C_1}\left(
\int_M R_{g'}^2e^{-\varphi}dV_{g'}\right).$$
Take $\varphi\in C^\infty(M)$ such that, for $g'=e^{-2\varphi}g$, $Vol(M,g')=1$ and such that
$\vert \nabla_g \varphi \vert_g \le C_1$ where $C_1$ is given by Proposition \ref{GE}.
Since $Vol(M,g')=1$, if $p$ is a point where $\varphi$ attains its minimum we have
$$e^{-3\varphi(p)} Vol(M,g) \ge 1,$$
and then, there exists $C_0$ depending only on $(M,g)$ such that
$\varphi(p) \le C_0$. Now, using the mean value theorem, it follows since
$\vert \nabla_g \varphi \vert_g$ is controlled by a constant depending only on
$(M,g)$, that $\max \varphi \le C_0'$ where $C'_0$ depends only on $(M,g)$.

\noindent Using this, we clearly have that
$$\int_M R_{g'}^2e^{-\varphi}dV_{g'} \ge e^{-C_0'} \int_M R_{g'}^2dV_{g'}.$$
Using H\"older inequality and the definition of the Yamabe invariant, we get (recall that
$Vol (M,g')=1$)
$$\int_M R_{g'}^2e^{-\varphi}dV_{g'} \ge e^{-C_0'}\left( Y(M,[g])\right)^2,$$
and then $I(M,g) \ge e^{-C_0'}\left( Y(M,[g])\right)^2$. This ends the proof.

\section{Lower bound}

\noindent For the lower bound, we need the
following lemmas:
\begin{lem}\label{change} For a conformal metric $\gt=e^{-2u}g$, we have the following integral
transformation
\bea
\int_M
\sdtu e^{-4u}\,dV_g&=&\int_M \sdu\,dV_g+\f{1}{8}\int_M R_g|\nabla_g
u|_g^2\,dV_g-\f{1}{4}\int_M|\nabla_g u|_g^4\,dV_g\\&&+\f{1}{2}\int_M\Delta_g
u|\nabla_g u|_g^2\,dV_g-\f{1}{2}\int_M A_g^1(\nabla_g u,\nabla_g u)\,dV_g.
\eea
\end{lem}
{\bf Proof} Denote $\sigut=\sutu$, $\sigu=\suu$, $\sigdt=\sdtu$,
$\sigd=\sdu$. We have
$$2\sigdt=\sigut^2-|A^1_{\gt}|^2_{\gt}.$$
By
equation (\ref{trasf2}), we have
$$\sigut e^{-2u}=\sigu+\Delta_g
u-\f{1}{2}|\nabla_g u|_g^2,$$
so
$$\sigut^2 e^{-4u}=\sigu^2+(\Delta_g u)^2
+\f{1}{4}|\nabla_g u|_g^4+2\sigu\Delta_g u-\Delta_g u|\nabla_g
u|_g^2-\sigu|\nabla_g u|_g^2.$$
After an easy computation, we get
\bea
|A^1_{\gt}|^2_{\gt}\,\,e^{-4u}&=&|A^1_g|_g^2+|\nabla_g^2
u|_g^2+\f{3}{4}|\nabla_g u|_g^4-\sigu|\nabla_g u|_g^2-\Delta_g u|\nabla_g u|_g^2+\\
&&+\,2(A^1_g)_{ij}\nabla_g^{2\, {ij}} u+2(A^1_g)_{ij}\nabla_g^i u\nabla_g^j
u+2\nabla_{g \, ij}^2 u\nabla_g^i u\nabla_g^j u.
\eea
Putting all together,
we obtain
\bea2
\sigdt e^{-4u}&=&2\sigd+(\Delta_g u)^2-|\nabla_g^2
u|_g^2-\f{1}{2}|\nabla_g u|_g^4+2\sigu\Delta_g u\\
&&-\,2(A^1_g)_{ij}\nabla_g^{2\, ij}u-2(A^1_g)_{ij}\nabla_g^i u\nabla_g^j
u-2\nabla^2_{g\, ij} u\nabla_g^i u\nabla_g^j u.
\eea
Now, by simple computation, we have the following identities
$$-2\int_M
(A^1_g)_{ij}\nabla_g^{2\, ij}u\,dV_g=-2\int_M\sigu\Delta_g u\,dV_g,$$
$$-2\int_M\nabla^2_{ij} u\nabla_g^i u\nabla_g^j u\,dV_g=\int_M\Delta_g u|\nabla_g u|_g^2\,dV_g,$$
where we integrated by parts and we used the Schur's Lemma for the
first identity. Finally we get
\bea2
\int_M\sigdt
e^{-4u}\,dV_g=2\int_M\sigd\,dV_g+\int_M\left[(\Delta_g u)^2-|\nabla_g^2 u|_g^2
-\f{1}{2}|\nabla_g u|_g^4 + \Delta_g u|\nabla_g u|_g^2-2A^1_g(\nabla_g u,\nabla_g
u)\right]\,dV_g,
\eea
Now using the integral Bochner formula
$$\int_M|\nabla_g^2 u|_g^2\,dV_g+\int_M Ric_g(\nabla_g u,\nabla_g
u)\,dV_g-\int_M(\Delta_g u)^2\,dV_g=0,$$
we get the final result.

\

\noindent In the sequel of the proof, we will need the following  proposition
(see \cite{GV} for the proof)
\begin{pro}\label{GV} If for some metric $g_1$ on $M$ we have
$A_{g_1}^t\in\Gamma^+_2$, then
\bea
-A_{g_1}^t+\sigma_1({g_1^{-1}}A_{g_1}^t)g_1&>&0,\\
A_{g_1}^t+\f{1}{3}\sigma_1({g_1^{-1}}A_{g_1}^t)g_1&>&0.
\eea
\end{pro}

\

\noindent Going on with the proof for the lower bound, we have the Lemma:

\begin{lem}\label{ineq}If $A_{\gt}^t\in\Gamma^+_2$, then we have the following estimate
$$\f{1}{2}\int_M
A_g(\nabla_g u,\nabla_g u)\,dV_g<\f{3-2t}{8}\int_M R_{\gt}|\nabla_g u|_g^2
e^{-2u}\,dV_g+\f{1}{4}\int_M\Delta_g u|\nabla_g
u|_g^2\,dV_g-\f{1}{4}\int_M|\nabla_g u|_g^4\,dV_g.$$
\end{lem}
{\bf Proof} Since $A^t_{\gt}\in\Gamma^+_2$, by Proposition
\ref{GV}, we get
$$-A^t_{\gt}>-\sigma_1(\gt^{-1}A^t_{\gt})\gt=-(4-3t)\sigma_1(\gt^{-1}A^1_{\gt})e^{-2u}g.$$
Hence we get
$$-A^1_{\gt}-(1-t)\sigma_1(\gt^{-1}A^1_{\gt})e^{-2u}g>-(4-3t)\sigma_1(\gt^{-1}A^1_{\gt})e^{-2u}g,$$
which implies that
$$A^1_{\gt}<(3-2t)\sigma_1(\gt^{-1}A^1_{\gt})e^{-2u}g.$$
Applying this to $\nabla_g u$ we obtain
$$\f{1}{2}A^1_{\gt}(\nabla_g u,\nabla_g u)<\f{3-2t}{8}R_{\gt}|\nabla_g u|_g^2e^{-2u}.$$
Using the conformal transformation law of the tensor $A_{\tilde g}^1$, integrating
over $M$, we have the result.

\

\noindent Now we are able to prove the following lower bound (recall
that $C_1$ is given by Lemma \ref{GE})

\begin{pro}[Lower Bound]\label{LB}
Assume that for some
$t\in[\delta,2/3]$   the following
estimate holds
\be
\int_M\sigd(g^{-1}A^1_g)\,dV_g+\f{1}{24}({7\over 10}-t)
\inf_{g'=e^{-2\varphi}g \, , \, \,  \vert \nabla_g \varphi\vert_g \le C_1}\left(
\int_M R_{g'}^2e^{-\varphi}dV_{g'}\right) := \mu_t>0.
\ee
Then there exists
$\underline{\delta}$ depending only on $(M,g)$ such that if $u_t\in
C^2(M)$ is a solution of (\ref{eq}) and if $A^t_{u_t}\in\Gamma^+_2$
then $u_t\ge\underline{\delta}$.
\end{pro}

\noindent{\bf Proof} Since $A_g^t=A_g^1+(1-t)\sigma_1(g^{-1}A_g^1)g,$ we easily have that
$$\sigma_2(A_g^t)=\sigma_2(A_g^1)+(1-t)(5-3t)\sigma_1(g^{-1}A_g^1)^2.$$
Letting
$\gt=e^{-2u_t}g$, \bea
e^{4u_t}f^2=\sigma_2(g^{-1}A^t_{u_t})&=&\sigma_2(g^{-1}A^1_{u_t})+(1-t)(5-3t)\left(\sigma_1(g^{-1}A^1_{u_t})\right)^2\\
&=&e^{-4u_t}\left(\sigma_2(\gt^{-1}A^1_{u_t})+\f{1}{16}(1-t)(5-3t)R_{\gt}^2\right).\eea
Integrating this with respect to $dV_g$ , we obtain
\bea
C\int_M
e^{4u_t}\,dV_g&\geq&\int_M f^2 e^{4u_t}\,dV_g\\
&=&\int_M
\sigma_2(\gt^{-1}A^1_{u_t})e^{-4u_t}\,dV_g+\f{1}{16}(1-t)(5-3t)\int_M
R_{\gt}^2 e^{-4u_t}\,dV_g\\
&=&\int_M
\sigma_2(\gt^{-1}A^1_{u_t})e^{-4u_t}\,dV_g+\f{1}{16}(1-t)(5-3t)\int_M
R_{\gt}^2 e^{-u_t}\,dV_{\gt},\eea where $C>0$ is chosen so that
$f^2\leq C$ (recall that, since $f=\sigma_2(g^{-1}A_g^\delta)$,
$C$ depends only on $(M,g)$). Using the fact that
$$R_{\gt}e^{-2u_t}=R_g+4\Delta_g
u_t-2|\nabla_g u_t|_g^2,$$
from Lemma \ref{change}, we get
\bea
\int_M
\sigma_2(\gt^{-1}A^1_{u_t})e^{-4u_t}\,dV_g&=&\int_M
\sigma_2(g^{-1}A^1_g)\,dV_g+\f{1}{8}\int_M R_{\gt}|\nabla_g
u_t|_g^2e^{-2u_t}\,dV_g\\&&-\f{1}{2}\int_M A_g^1(\nabla_g u,\nabla_g
u)\,dV_g.
\eea
Notice that, since  $A^t_{u_t} \in \Gamma_2^+$, we have
$$0<\sigma_1(g^{-1}A^t_{u_t})=(4-3t)\sigma_1(g^{-1}A^1_{u_t}),$$
and so $R_{\gt}>0$. By
Lemma \ref{ineq}, we obtain
\bea
\int_M
\sigma_2(\gt^{-1}A^1_{u_t})e^{-4u_t}\,dV_g&\geq&\int_M
\sigma_2(g^{-1}A^1_g)\,dV_g-\f{1-t}{4}\int_M R_{\gt}|\nabla_g
u_t|_g^2e^{-2u_t}\,dV_g\\&&-\f{1}{4}\int_M\Delta_g u_t|\nabla_g
u_t|_g^2\,dV_g+\f{1}{4}\int_M|\nabla_g u_t|_g^4\,dV_g.
\eea
By Young's
inequality, one has
$$\int_M R_{\gt}^2
e^{-u_t}\,dV_{\gt}\geq\f{2}{\varepsilon}\int_M R_{\gt}|\nabla_g
u_t|_g^2e^{-2u_t}\,dV_g-\f{1}{\varepsilon^2}\int_M|\nabla_g u_t|_g^4\,dV_g,$$
for all $\varepsilon>0$. By an easy computation, we have
$$\f{1}{16}(1-t)(5-3t)=\f{1}{24}({7\over 10}-t)+P_2(t),$$
where $P_2(t)$
is a positive, second order, polynomial in $t$. Putting all
together, we obtain (for $C>0$ depending only on $(M,g)$)
\bea
C\int_M e^{4u_t}\,dV_g&\geq&\int_M
\sigma_2(\gt^{-1}A^1_{u_t})e^{-4u_t}\,dV_g+\f{1}{16}(1-t)(5-3t)\int_M
R_{\gt}^2 e^{-u_t}\,dV_{\gt}\\&=&\int_M
\sigma_2(\gt^{-1}A^1_{u_t})e^{-4u_t}\,dV_g+\left(\f{1}{24}({7\over 10}-t)+P_2(t)\right)\int_M
R_{\gt}^2 e^{-u_t}\,dV_{\gt}\\&\geq&\int_M
\sigma_2(g^{-1}A^1_g)\,dV_g+\f{1}{24}({7\over 10}-t)\int_M R_{\gt}^2
e^{-u_t}\,dV_{\gt}\\&&+P_2(t)\int_M R_{\gt}^2
e^{-u_t}\,dV_{\gt}-\f{1-t}{4}\int_M R_{\gt}|\nabla_g
u_t|_g^2e^{-2u_t}\,dV_g\\&&-\f{1}{4}\int_M\Delta_g u_t|\nabla_g
u_t|_g^2\,dV_g+\f{1}{4}\int_M|\nabla_g u_t|_g^4\,dV_g.
\eea
Now using Young's
inequality  and the conformal change equation of the scalar
curvature, we get (for a certain $C>0$ depending only on $(M,g)$)
\bea
C\int_M e^{4u_t}\,dV_g&\geq&\int_M
\sigma_2(g^{-1}A^1_g)\,dV_g+\f{1}{24}({7\over 10}-t)\int_M R_{\gt}^2
e^{-u_t}\,dV_{\gt}\\&&+\left(\f{2P_2(t)}{\varepsilon}-\f{1-t}{4}\right)\int_M
R_g|\nabla_g
u_t|_g^2\,dV_g\\&&+\left(\f{8P_2(t)}{\varepsilon}-(1-t)-\f{1}{4}\right)\int_M\Delta_g
u_t|\nabla_g
u_t|_g^2\,dV_g\\&&+\left(\f{3-2t}{4}-\f{P_2(t)}{\varepsilon^2}-\f{4P_2(t)}{\varepsilon}\right)\int_M|\nabla_g
u_t|_g^4\,dV_g.
\eea
We choose $\varepsilon=\varepsilon(t)>0$, such that
$\f{8P_2(t)}{\varepsilon}-(1-t)-\f{1}{4}=0$. One can easily check that, with this choice,
$$\f{2P_2(t)}{\varepsilon}-\f{1-t}{4} \ge 0 \quad \hbox{ and } \quad
\f{3-2t}{4}-\f{P_2(t)}{\varepsilon^2}-\f{4P_2(t)}{\varepsilon} \ge 0.$$
Finally, recalling that according to lemma \ref{GE} $\Vert \nabla_g u_t \Vert_{g, \infty} \le C_1$
with $C_1$ depending only on $(M,g)$,
we obtain the following
estimate (for a certain $C>0$ depending only on $(M,g)$)
\bea
C\int_M e^{4u_t}\,dV_g
&\geq&\int_M
\sigma_2(g^{-1}A^1_g)\,dV_g+\f{1}{24}({7\over 10}-t)\int_M R_{\gt}^2
e^{-u_t}\,dV_{\gt} \\
&\ge&
\int_M\sigd(g^{-1}A^1_g)\,dV_g+\f{1}{24}({7\over 10}-t)
\inf_{g'=e^{-2\varphi}g \, , \, \,  \vert \nabla_g \varphi\vert_g \le C_1}\left(
\int_M R_{g'}^2e^{-\varphi}dV_{g'}\right)=\mu_t >0.
\eea
This gives
$$\max_M u_t\geq\log\mu_t-C(g).$$
Since $\Vert\nabla_g u_t\Vert_{g,\infty}<C_1$ this implies
the Harnack inequality
$$\max_M u_t\leq\min_M u_t+C(M,g),$$
by
simply integrating along a geodesic connecting points at which $u_t$
attains its maximum and minimum. Combining this two inequalities, we
obtain $$\min_M u_t\geq\log\mu_t-C,$$
where $C$ only depends on $(M,g)$. This ends the proof of the Lemma.

\smallskip
\section{$C^{2,\alpha}$ estimate}
We have the following  $C^{2,\alpha}$ estimate for
solutions of the equation (\ref{eq1}). For the proof, see \cite{GV}
and \cite{GV2}.

\begin{pro}[$C^{2,\alpha}$ estimate]\label{C2}
Let $u_t\in C^4(M)$ be a
solution of (\ref{eq}) for some $\delta\leq t\leq 2/3$, satisfying
$\underline{\delta}<u_t<\bar\delta$, and $\parallel\nabla
u_t\parallel_{{g,\infty}}<C_1$. Then for $0<\alpha<1$, $\parallel
u_t\parallel_{g,C^{2,\alpha}}\leq C_2$, where $C_2$ depends only on
$(M,g)$.
\end{pro}

\section{Proof of Theorem \ref{Tgen}}

\noindent We use the continuity method. Our 1-parameter family of equations,
for $t\in[\delta,t_0]$, is
\be\label{eqaz}\sigma^{1/2}_2(g^{-1}A^t_{u_t})=f(x)e^{2u_t},\ee with
$f(x)=\sigma^{1/2}_2(g^{-1}A^\delta_g)>0,$ and $\delta$ was chosen
so that $A^\delta_g$ is positive definite. Define
$$\mathcal{S}=\left\{t\in[\delta,t_0]\mid\exists\mbox{\,a
solution\,}u_t\in
C^{2,\alpha}(M)\mbox{\,of\,}(\ref{eqaz})\mbox{\,with\,}A^t_{u_t}\in\Gamma^+_2\right\}.$$
Clearly, with our choice of $f$, $u\equiv 0$ is a solution for
$t=\delta$. Since $A^{\delta}_g$ is positive definite,
$\delta\in\mathcal{S}$, and $\mathcal{S}\neq\emptyset.$ Let
$t\in\mathcal{S}$, and $u_t$ be a solution. By Proposition
\ref{Inver}, the linearized operator at $u_t$,
$\mathcal{L}^t:C^{2,\alpha}(M)\rightarrow C^{\alpha}(M)$, is
invertible. The implicit function theorem  tells us
that $\mathcal{S}$ is open. From classical elliptic theory, it
follows that $u_t\in C^{\infty}(M)$, since $f\in C^{\infty}(M)$. By
Proposition \ref{UB} we get an uniform upper bound on the
solutions $u_t$, independent of $t$. We may then apply Proposition
\ref{GE} to obtain a uniform gradient bound on $u_t$, and by
Proposition \ref{LB}, we get a uniform lower bound. Finally using
Proposition \ref{C2} and the classical Ascoli-Arzela's Theorem,
then implies that $\mathcal{S}$ must be closed, therefore
$\mathcal{S}=[\delta, t_0]$. The metric $\gt=e^{-2u_{t_0}}g$ then
satisfies $\sigma_2(A^{t_0}_{\gt})>0$ and $R_{\gt}>0$. The
inequalities \eqref{inequal} follow from proposition \ref{GV}.

\section{Proof of Theorem \ref{T2}}

\noindent Theorem \ref{T2} is a direct consequence of Theorem \ref{Tgen} and of
Lemma \ref{techn}.

\section{Proof of Corollary \ref{Qcurv}}

Assume that $M$ admits a metric $g'$ such that
$Q_{g'}\ge {1\over 48} R^2_{g'}$ and $Y(M,[g'])\ge 0$. Recall that
$$Q_{g'} = -{1\over 4}\Delta_{g'} R_{g'} - 2 \vert Ric_{g'} \vert_{g'}^2 +{23\over 32}R_{g'}^2,$$
Integrating $Q_{g'}$ on $M$ with respect to $dV_{g'}$ we obtain
(since $Q_{g'}\ge 0$)
\be\label{control}
\int_M \vert Rig_{g'} \vert_{g'}^2 dV_{g'} \le {23 \over 64} \int_M R_{g'}^2 dV_{g'}.
\ee
Now if we compute $\int_M \sigma_2(g'^{-1}A^1_{g'})$ using \eqref{control},
 we have (recall that $\sigma_2(g'^{-1}A^1_{g'})=-{1\over 2}\vert Ric_{g'} \vert^2_{g'}
+{3\over 16} R^2_{g'}$):
$$\int_M \sigma_2(g'^{-1}A^1_{g'}) \ge {1\over 128}\int_M R^2_{g'}dV_{g'} \ge 0.$$
Now, consider the conformal laplacian operator
$L_{g'} := \Delta_{g'} -{1\over 8} R_{g'}$. We have using the assumption
 $Q_{g'} \ge {1\over 48}R^2_{g'}$
$$L_{g'}R_{g'} = \Delta_{g'}R_{g'} -{1\over 8} R^2_{g'} \le
-8\vert Ric_{g'} \vert^2_{g'} + {22\over 8}R^2_{g'} -{1\over 12}R^2_{g'}\le
\left(-{8\over 3} +{22\over 8}-{1\over 12}\right)R^2_{g'} =0.
$$
Applying a Lemma due to Gursky \cite{Gurann}, since
$Y(M,[g'])\ge0$ we have either $R_{g'}>0$ (if $Y(M,[g'])>0$) or
$R_{g'}\equiv 0$ (if $Y(M,[g'])=0$).
If $Y(M,[g'])>0$ we can apply Theorem \ref{T1}
to conclude that $m$ is diffeomorphic to a spherical space form. Otherwise, if $Y(M,[g'])=0$, since
$Q_{g'}\ge {1\over 48} R^2_{g'}$ and $R_{g'}\equiv 0$, we deduce, using the
expression giving $Q_{g'}$, that $Ric_{g'}\equiv 0$ and then $M$ is diffeomorphic to a
quotient of ${\Bbb R}^3$.

\noindent This ends the proof of the Corollary.

\end{document}